\documentclass[12pt,psamsfonts]{article}

\usepackage[margin=1in]{geometry}
\usepackage{amsmath,amsfonts,amssymb,amsmath,amscd}
\usepackage[active]{srcltx}
\usepackage{graphicx,array}
\usepackage{amsthm}
\usepackage{amsbsy}
\usepackage{tikz}
\usepackage{float}

\newcolumntype{L}{>{\displaystyle}l}
\newcolumntype{C}{>{\displaystyle}c}
\newcolumntype{R}{>{\displaystyle}r}

\numberwithin{equation}{section}
\numberwithin{figure}{section}

\newtheorem{theorem}{Theorem}[section]
\newtheorem{lemma}[theorem]{Lemma}
\newtheorem{proposition}[theorem]{Proposition}
\newtheorem{corollary}[theorem]{Corollary}
\newtheorem{definition}[theorem]{Definition}
\newtheorem{example}[theorem]{Example}
\newtheorem{remark}[theorem]{Remark}

\newcommand{\dem}{\noindent{\rm Proof:\ }}

\newcommand{\eproof}{\noindent\mbox{\framebox [0.6ex]{}}  \medskip}

\def \r{\mbox{${\mathbb R}$}} 
\def \z{\mbox{${\mathbb Z}$}} 
\def \n{\mbox{${\mathbb N}$}} 
\def \S{\mbox{\rm{${\textbf{S}}$}}}
\def \R{\mbox{\rm{${\textbf{R}}$}}}
\def \T{\mbox{\rm{${\textbf{T}}$}}}
\def \Z{\mbox{\rm{${\textbf{Z}}$}}}

\def \G{\mbox{${\mathcal{G}}$}}

\def \E{\mbox{${\mathcal{E}}$}}
\def \V{\mbox{${\mathcal{V}}$}}
\def\cN{\mathcal{N}}

\begin{document}

\title{Gradient systems on coupled cell networks}

\author{Miriam Manoel\\
{\small Department of Mathematics, ICMC - University of S\~ao Paulo}\\
{\small C.P. 668, 13560-970 S\~ao Carlos SP, Brazil \footnote{Email address: miriam@icmc.usp.br (corresponding author)}}
\and Mark Roberts\\
{\small Department of Mathematics, University of Surrey}\\
{\small Guildford, GU2 7XH, UK \footnote{Email address: m.roberts@surrey.ac.uk \hfill  30 January 2015}}\\
{\small and}\\
{\small African Institute of Mathematical Sciences, Tanzania}\\
{\small PO Box 1200, Arusha, Tanzania}} 
\date{}

\maketitle

\begin{flushleft}

\vspace{2mm}

\end{flushleft}

\begin{abstract}
For networks of coupled dynamical systems we characterize admissible functions, that is, functions whose gradient is an admissible vector field. The schematic representation of a gradient network dynamical system is of an undirected cell graph, and we use tools from graph theory to deduce the general form of such functions, relating it to the topological structure of the graph defining the network. The coupling of pairs of dynamical systems cells is represented by edges of the graph, and from spectral graph theory we detect the existence and nature of equilibria of the gradient system from the critical points of the coupling function. In particular, we study fully synchronous and 2-state patterns of equilibria on regular graphs.These are two special types of equilibrium configurations for gradient networks. We also investigate equilibrium configurations of ${\bf S}^1$-invariant admissible functions on a ring of cells. 
\end{abstract}


\section{Introduction} \label{sec: introduction}

Networks of coupled dynamical systems frequently demonstrate phenomena, such as synchrony, 
phase relations, resonances and non-generic bifurcations, that are typically associated with group equivariance, even when the networks are not invariant 
under any non-trivial group actions. They appear to have some form of `hidden symmetry'. These phenomena have been explored from a number of different points of view, {\it e.g.} using groupoid formalism, as in \cite{Field} and \cite{Gol Stewart} or, more generally, category theory, as in \cite{De Ville Lerman}, and semigroup formalism, as in \cite{Rink Sanders}.  The algebraic formalism of symmetry groupoids of networks has led to a successful way to search for patterns of synchrony  on networks. In general, elements of the groupoid  can be thought of as a set of local symmetries of the network, relating sets of cells in the network
in a way that the system is invariant under the action of these symmetries. This has been first established in \cite{SGP} and explored since in a great number of works; for example,  \cite{AS, DS, WG},   among many others.  

In this paper we focus on one particular class of systems, 
namely gradient systems. Synchrony manifests itself in a network as configurations of the coupled systems that behave setwise identically. For gradient systems there are also configurations given by phase relations between cells, which  appear naturally in a class of vector fields generated by functions with an extra ${\bf S}^1$-invariance. 

We consider a network of a finite set of dynamical systems cells coupled together in a manner given schematically by a connected undirected graph $\G =(\V, \E)$ without multiple edges, whose vertices represent the cells and edges correspond to couplings. We shall call $\G$ a cell graph. The set $\V = \{v_1,\ldots, v_n\}$ denotes the set of vertices and $\E$ the set of edges of $\G$, which  are identified as 2-element subsets of $\V$ although we shall denote an element of $\E$ as a pair $(v_1, v_2)$. We assume that a cell graph may contain all internal edges (or loops), in which case  $\V \subseteq \E$. We shall also denote by $I(v)$ the {\it input set} of $v \in \V$, that is, the set of
vertices $u \in \V$, $u \neq v$, such that $(u,v) \in \E$, and $d(v)$ shall denote the {\it degree} of $v$, the number of vertices in $I(v)$. Notice that for our purposes we shall not consider  loops to be in the input sets.  A smooth 
manifold $P_v$ is assigned to each vertex $v \in \V$, so that the total configuration space is
$P =  \prod_{v \in \V} P_v.$  
We denote this network by $\cN = (\G, P).$ 

A vector field $g = (g_1, \ldots, g_n)$ on $P$ is called an admissible vector field  if it is consistent with the network structure, that is, if it is equivariant under the action of  the network symmetry groupoid  (Definition 4.1 in \cite{SGP}). In particular, if $x_v$ denotes coordinates on $P_v$, each component $g_{v}$ must depend only on the variables $x_u$ for which $(u, v) \in \E$. More precisely, for $ I(v) \ = \ \{v_1, \ldots, v_{d(v)}\},$
\begin{equation} \label{eq: dependency}
 g_v (x) = \tilde{g}_v(x_v, x_{v_1}, \ldots, x_{v_{d(v)}}). 
\end{equation} 
In addition, equality constraints are imposed 
between components corresponding to cells $u$ and $v$  whose input sets $I(u)$ and $I(v)$ are isomorphic by an element of the groupoid, that is, 
\begin{equation} \label{eq: equality constraint}
\tilde{g}_u =  \tilde{g}_v.
\end{equation} 
General network dynamical systems are usually defined for  directed graphs. For a network of smooth dynamical systems to be given by the  (negative) gradient  of 
a smooth function $f: P \to \r$,
\begin{equation} \label{eq: gradient}
\dot{x} \ = \ - \nabla f(x), 
\end{equation}
it is necessary that if $(u,v)$ is an edge then so is $(v,u)$, and hence the directed graph is equivalent to an undirected graph.
We assume that all edges are identical and all vertices are of the same type. In particular, all components $P_v$ of $P$ are the same. Since $\nabla f$ is an admissible vector field, then its components satisfy (\ref{eq: dependency}) and
(\ref{eq: equality constraint}).  Since we are assuming that all cells are identical, then (\ref{eq: equality constraint}) holds precisely if 
$ d(u) = d(v)$, and, for each $v \in \V,$
\begin{equation} \label{eq: first partial derivatives}
 \frac{\partial f}{\partial x_v} (x) \ =  \ \tilde{f}_{d(v)}(x_v, \overline{x_{u_1}, \ldots, x_{u_{d(v)}}}), 
\end{equation}
where the over-bar  indicates  invariance by permutation of these variables. 

Simple examples show that not all admissible vector fields are gradient. In this paper we restrict to the class of admissible vector fields that are of gradient type. Our main result, Theorem~\ref{thm: admissible functions}, characterizes the admissible functions, namely, smooth functions whose gradients are admissible vector fields. As we shall see, for any network these functions are decomposed as a sum of  components that depend on each cell individually - the self-connection functions - and those that depend on the way the cells are coupled together - the coupling functions.  

This work also addresses the analysis of critical points of admissible functions $f$ or, equivalently, the analysis of equilibria of (\ref{eq: gradient}).   
As suggested by the general form of these functions, there is a direct relationship between critical points of the associated coupling function
and critical points of the admissible function that are either totally synchronous (all cells assume the same value) or given by 2-colour patterns (when each cell assumes one out of two possible distinct values). 
In fact, this makes these two configurations on networks special in the class of gradient systems.  In this paper we investigate these two particular types of critical points, regarding existence and nature, and shall understand how these are related to both the critical points of the coupling functions and the architecture of the graph of couplings.   This study does not generalize for patterns with more than two colours. For these cases a different general  approach can be applied, and this is done in \cite{ADM}. In the presence of extra symmetries, however, other types of critical points of admissible functions become expected for some classes of graphs. These can appear with a variety of configurations, with  no nontrivial synchrony among cells for example, but still with a  phase relation. Specifically, this is the case when the admissible function is invariant under the circle group ${\bf S}^1.$

There are similarities between the network gradient dynamical systems that we consider and Ising-Potts model, Kuramoto model, antiferromagnetic XY model (also called AFXY model), neural networks  and other interaction systems. For example, applications of our results given in Section~\ref{sec:S1-invariant admissible functions} are related to
the results of  \cite{Ermentrout}, where the author applies the method of averaging to reduce models of discrete and a continuum of  neural arrays to systems of phase equations, that is, to equations in which each cell is represented by a single
variable lying on ${\bf S}^1$ and interactions between two connected neural cells  are periodic functions that depend only on the difference between their two phases. In such model, an edge  $(u,v) \in \E$ represents a synapse between the two neurons $u$ and $v.$ There is also a direct connection between our study of critical points of admissible functions given in Subsection~\ref{subseq: minimum on a ring} and several results about ground states in Kuramoto and AFXY models (see  \cite{Bronski et al, Korshunov, Lee, MS}).  

This paper is organized as follows: In Section~\ref{sec: admissible functions} we present the general form of admissible functions on networks, which is our main result. In Section~\ref{sec:reg graphs}  we study critical points of admissible functions, relating their existence and nature with the corresponding critical points of the coupling functions. A detailed analysis of fully synchronous critical points and 2-colour patterns of critical points is carried out for the special case of regular graphs. Section~\ref{sec:S1-invariant admissible functions} presents a study of all possible critical configurations of a class of admissible functions under an extra ${\bf S}^1$-invariance.


\section{Admissible functions} \label{sec: admissible functions}

For a given cell graph $\G = (\V, \E)$,  we characterize the functions 
on $P =  \prod_{v \in \V} P_v$ whose gradients are admissible vector field for $\cN=(\G,P).$   We shall consider 
$P_v = {\r}^k$, $k \geq 1$, for all $v \in {\V}$, taking local coordinates if necessary.  The variable on each cell $v \in \V$ shall be denoted by $x_v = (x_v^1, \ldots, x_v^k)$. 

\begin{definition} \label{def: admissible function}
For a cell network $\cN = (\G, P)$, a smooth function 
$f: P \to \r$  is an admissible function if its gradient $\nabla f$ is an admissible vector field for $\cN$.
\end{definition}

The next result is one of the main constraints imposed on admissible vector fields if they are of gradient type. 

\begin{lemma} \label{lem: third derivatives}
For a cell network $\cN = (\G, P)$ with $|\V| \geq 3$, if  a smooth function 
$f: P \to \r$ is an admissible function, then, unless coupling is all-to-all, we have, for any
$1 \leq i,j,l \leq k$, 
\begin{equation} \label{eq: third partial derivatives}
 \frac{\partial^3 f}{\partial x_{v_1}^i \partial x_{v_2}^j \partial x_{v_3}^l} \equiv 0,
\end{equation}
if $v_1,v_2,v_3$ are distinct vertices in $\G$.
\end{lemma}
\dem 
We prove by induction on the number of vertices in $\G$. For simplicity the proof is carried out  for $k = 1$, but the result holds equally for any higher dimension $k>1$.  For $|\V|=3$,  if (\ref{eq: third partial derivatives}) does not hold, then for each $i$, $i=1,2,3$, $\partial f \slash \partial x_{v_i}$ is a nontrivial 
function of the other two variables, and so $\G$ necessarily has its three vertices coupled one another. 
For $n \geq 4$, assume that if a network has $n-1$ cells and three distinct cells  for which   (\ref{eq: third partial derivatives}) does not hold, then these are all coupled one another. Now,  let  $|\V|= n$ and 
suppose that there exist distinct $v_1, v_2, v_3 \in \V$ such that  (\ref{eq: third partial derivatives})  fails.
Let $u \in \V$, $u  \neq v_1, v_2, v_3$.
Since $\G$ is connected, there exists $w \in \V$ such that
$ (u,w) \in \E$. By the induction assumption all vertices in $\G$, except possibly $u,$ are all-to-all coupled. In particular $w \in I(v_i)$, $i=1, 2, 3$.  Now, for $i=1, 2, 3$, we use the form (\ref{eq: first partial derivatives}) of $\partial f \slash \partial x_{v_i}$ together with the hypothesis of nonvanishing third-order derivatives with respect to $x_{v_1}, x_{v_2}, x_{v_3}$  to conclude that 
\begin{equation} \label{eq: inequality}
 \frac{{\partial}^{3} f}{\partial x_{v_1} \partial x_{v_2} \partial x_{w}} \not\equiv 0. 
\end{equation}
But by  the induction hypothesis we also have
\[ \frac{\partial f}{\partial x_w} \ = \ \tilde{f}_{n-1}(x_w, \overline{x_{v_1}, \ldots, x_{v_{n-2}},x_u}). \]
So for all $i=1, \ldots, n-2$, 
\[ \frac{\partial^{2}}{\partial x_{u} \partial x_{w}} \bigl(\frac{\partial f}{ \partial x_{v_i}}\bigr) \not\equiv 0, \]
and hence 
$ u \in  \cap_{v \in \V} I(v), $
that is, vertices in $\G$ are all-to-all coupled. If $u= v_i$, for some $i \in \{1,2,3\}$, we just rearrange
the indices in (\ref{eq: inequality}) to get the same conclusion.  \hfill \eproof

We set up the following notation: for each edge $e = (v_1,v_2) \in \E$ representing the coupling of vertices $v_1, v_2 \in \V$, we
denote $\rho(e), \tau(e) \in \{v_1, v_2\},$ where $\rho(e) \neq \tau(e)$ if $v_1 \neq v_2$ and write the ordered pair $(\rho(e), \tau(e) )$ to represent a directed edge, with `head' $\rho(e)$ and `tail' $\tau(e),$ corresponding to the (undirected) edge $e$.  

As a consequence of Lemma~\ref{lem: third derivatives}, an admissible function is of the form
\begin{equation} \label{eq: double sum1}
f(x) \ = \ \sum_{e \in \E} \beta_e(x_{\rho(e)}, x_{\tau(e)}) \ + \ \sum_{v \in \V} \alpha_v(v),
\end{equation}
for smooth functions $\beta_{e}: \r^{2k} \to \r$ and $\alpha_v: \r^k \to \r$. We assume without loss of generality that the functions $\beta_e$'s do not contain terms of just one variable, and that $f$ vanishes at the origin. For $e \in \E$ and $v \in \V$, the functions $\beta_e$'s   in the first sum shall be called {\it coupling} functions, and the $\alpha_v$'s shall be called {\it self-connection} functions. 

Theorem~\ref{thm: admissible functions} below characterizes admissible functions defined on networks of coupled cells. As we shall see, these are distinguished when $\G$ is a bipartite graph. 

\begin{remark} \label{rmk: bipartite graph}
{\rm Recall that a graph $\G$ is bipartite if its set of vertices $\V$ can be devided into two disjoint subsets $\V_1$ and $\V_2$ such that every edge of $\G$ connects a vertex in a subset to a vertex in the other. In this case,  a 2-colouring can be defined on  the graph (\cite{ADH}) such that  any two vertices in either $\V_1$ or $\V_2$ receive the same colour.  For graphs of many cells, algorithms can be useful to check whether a graph is or is not bipartite (see \cite{Kleinberg} and \cite{Sedgewick}).}\end{remark}

For any subgraph ${\cal S}$ of $\G$ let us denote $\E_{\cal S}$ its set of edges. Also, let ${\bf S}_n$ denote the permutation group acting on $(\R^k)^n$ by permutation of variables and  let ${\bf Z}_2$ denote the order-2 permutation group acting on any subspace $(\r^{k})^2$ of two variables in $(\R^k)^n$.

\begin{theorem} \label{thm: admissible functions} 
If $\cN =(\G, P)$ is an all-to-all coupling network of $n$ cells, then admissible functions are the $\S_n$-invariant functions.
Otherwise, a function  $f: P \to \r$ is an  admissible function associated to $\cN$ if, and only if, 
there exist smooth functions $\alpha_{d(v)} : \r^k \to \r$ and
$\beta : \r^{2k}  \to \r$ such that one of the following holds: 
\begin{enumerate} 
\item If $\G$ is bipartite, then for the disjoint union  $\V = \V_1 \mathbin{\dot{\cup}} \V_2$ we have
\begin{equation} \label{eq: eq1 thm admissible functions}
 f(x) \ = \   \sum_{e \in \E, \ \rho(e) \in \V_i} \beta(x_{\rho(e)}, x_{\tau(e)}) + \sum_{v \in \V_1} \alpha_{d(v)}(x_{v}) + \sum_{v \in \V_2} \gamma_{d(v)}(x_{v}), 
\end{equation}
where $i=1$ or $2$. In addition, if there exist $v_1 \in \V_1$ and $v_2 \in \V_2$ such that
$d(v_1) = d(v_2)$, then $\beta$ is  ${\bf Z}_2$-invariant and
\begin{equation} \label{eq: eq2 thm admissible functions}
 f(x) \ = \  \sum_{e \in \E} \beta(x_{\rho(e)}, x_{\tau(e)}) + \sum_{v \in \V} \alpha_{d(v)}(x_{v}). 
\end{equation}
\item If $\G$ is not bipartite, then $f$ is of the form (\ref{eq: eq2 thm admissible functions}) with
$\beta$ ${\bf Z}_2$-invariant.
\end{enumerate}
\end{theorem}
\dem If the cells in $\G$ are all coupled one another, then $\nabla f$ is $\S_n$-equivariant, so $f$ is
$\S_n$-invariant. 

It is easy to verify that functions of the forms (\ref{eq: eq1 thm admissible functions}) and
(\ref{eq: eq2 thm admissible functions}) are admissible functions. For the converse,  since we are interested in the non all-to-all coupling case,  we have $n \geq 3$. We present the proof for $k=1$, since for $k>1$ the proof adapts straightforwardly. In fact, the arguments rely essentially on results from graph theory and the repeated use of (\ref{eq: first partial derivatives}); if $k > 1$, then for $x_v = (x_{v}^1, \ldots, x_v^k),$ (\ref{eq: first partial derivatives}) is 
\[ \frac{\partial f}{\partial x_v^i} (x) \ =  \ \tilde{f}^i_{d(v)}(x_v, \overline{x_{u_1}, \ldots, x_{u_{d(v)}}}), \ \ i=1, \ldots, k, \]
and we work out the partial derivatives on each component of $x_v$ to obtain the same conclusions as for $k=1$. 

Assume  that $\G$ is bipartite.  Considering the partition $\V = \V_1 \mathbin{\dot{\cup}} \V_2,$ we have a decomposition 
\[ \G = {\mathbin{\dot{\cup}}}_{v \in \V_1} {\cal K}^v, \]
where ${\cal K}^v$ denotes the star graph $K_{1, d(v)}$ with centre vertex $v \in \V_1$ and $d(v)$ edges $(v,u)$, for all $u \in I(v)$. We notice that the disjoint union above refers to non-repeating edges, and  it can obviously be taken to run through elements in either subset of the partition, indistinctly. Using the decomposition above, we rewrite (\ref{eq: double sum1}) as
\[ f(x)  =   \sum_{v \in \V_1} \sum_{e \in \E_{{\cal K}^v}} \beta_e(x_{\rho(e)}, x_{\tau(e)}) \ 
+ \ \sum_{v \in \V} \alpha_v(x_v) \ =  \sum_{v \in \V_1} \sum_{u \in I(v)} \beta_{uv}(x_v, x_u) \ + \
 \sum_{v \in \V} \alpha_v(x_v).  \]
From (\ref{eq: first partial derivatives}), for each $v \in \V_1$ with $d(v) \geq 2$, the partial derivative $\partial f/ \partial x_v$ is invariant by pairwise permutation of $u_i, u_j \in I(v).$ Hence,
\[ \frac{\partial \beta_{u_iv}}{\partial x_v} (x_v, x_{u_1}) \ = \ \frac{\partial \beta_{u_jv}}{\partial x_v} (x_v, x_{u_1}). \]
Since coupling functions do not depend on each variable independently, the equality above yields  $\beta_{u_iv} = \beta_{u_jv}$. Thus,
\[ f(x) \ = \  \sum_{v \in \V_1} \sum_{u \in I(v)} \beta_{v}(x_v, x_u) \ + \
 \sum_{v \in \V} \alpha_v(x_v).  \]
If $\G$ is the star graph, then $\V_1= \{v\}$. In this case, by (\ref{eq: first partial derivatives}),
\[ \frac{\partial \beta_v}{\partial y}(x_v, x_{u_i}) + \alpha^{\prime}_{u_i}(x_{u_i}) \ = \  
\frac{\partial \beta_v}{\partial y}(x_v, x_{u_j}) + \alpha^{\prime}_{u_j}(x_{u_j}), \]
for any  pair $u_i, u_j \in I(v)$. Set $x_{u_i} = x_{u_j}$ to conclude that $\alpha_{u_i} = \alpha_{u_j}.$ Hence, the self-connection functions in $\V_2$ are all equal, and so in this case $f(x) = \sum_{e \in \E} \delta(x_{\rho(e)}, x_{\tau(e)})$, for some smooth function $\delta$.  

If $\G$ is not the star graph, then there exists $u \in \V_2$ with $d(u) \geq 2$. Consider then two edges  $(u, v_1), (u, v_2) \in \E$. We have $\partial f/ \partial x_u$  invariant by interchanging $v_1$ and $v_2$, so
\[ \frac{\partial \beta_{v_1}}{\partial y} (x_{v_1}, x_u) +  \frac{\partial \beta_{v_2}}{\partial y} (x_{v_2}, x_u) = \frac{\partial \beta_{v_1}}{\partial y} (x_{v_2}, x_u) +  \frac{\partial \beta_{v_2}}{\partial y} (x_{v_1}, x_u), \] 
which yields $\beta_{v_1} = \beta_{v_2}$. Since $\G$ is connected, we apply this idea transitively through all edges of $\G$ to conclude that
\[ f(x) \ = \ \sum_{v \in \V_1} \sum_{ u \in I(v)} \beta(x_v, x_u) + \sum_{v \in \V} \alpha_v(x_v). \]
Now let $v_1, v_2 \in \V_1$ with $d(v_1) = d(v_2)$. For any $a, b$, take 
\begin{equation} \label{eq: x bar}
\bar{x} = (x_{v_1}, \ldots, x_{v_n}), \  \ x_{v_1} = x_{v_2} = a, \ \ u = b, \ 
\forall u \in I(v_1) \cup I(v_2). 
\end{equation}
By (\ref{eq: first partial derivatives}) applied to $\partial f / \partial x_{v_1}$ and
$\partial f / \partial x_{v_2}$ at $\bar{x}$, we have
\[ d(v_1) \frac{\partial \beta}{\partial x}(a,b) + \alpha^{\prime}_{v_1}(a) =   
d(v_2) \frac{\partial \beta}{\partial x}(a,b) + \alpha^{\prime}_{v_2}(a), \]
so $\alpha_{v_1} = \alpha_{v_2}$. Analogously, for $v_1, v_2 \in \V_2$ with
$d(v_1) = d(v_2)$ we obtain
\[ d(v_1) \frac{\partial \beta}{\partial y}(b,a) + \alpha^{\prime}_{v_1}(a) =   
d(v_2) \frac{\partial \beta}{\partial y}(b,a) + \alpha^{\prime}_{v_2}(a), \]
and so $\alpha_{v_1} = \alpha_{v_2}$. Therefore, self-coupling functions are invariant by degree on each subset of the partition of $\V$, yielding  (\ref{eq: eq1 thm admissible functions}).

Suppose now that there exist $v_1 \in \V_1$ and $v_2 \in \V_2$ with  $d(v_1) = d(v_2) = d$. Apply (\ref{eq: first partial derivatives})  to $\partial f / \partial x_{v_1}$ and
$\partial f / \partial x_{v_2}$ at $\bar{x}$ given in (\ref{eq: x bar}) to obtain
\[ d \  \frac{\partial \beta}{\partial x}(a,b) + \alpha^{\prime}_{v_1}(a) =   
d  \frac{\partial \beta}{\partial y}(b,a) + \alpha^{\prime}_{v_2}(a), \]
that is,
\[ d \bigl(  \frac{\partial \beta}{\partial x}(a,b) - \frac{\partial \beta}{\partial y}(b,a) \bigr) \ = \ {(\alpha_{v_2} - \alpha_{v_1})}^{\prime} (a). \]
Since $\beta$ does not depend on each variable independently, the above equality implies that $\beta$ is ${\bf Z}_2$-invariant and $\alpha_{v_1} = \alpha_{v_2}$, resulting in  (\ref{eq: eq2 thm admissible functions}). 

Assume now that  $\G$ is not bipartite. We consider a spanning tree ${\cal T}$ of 
$\G$ to rewrite (\ref{eq: double sum1}) as
\[ f(x) \ = \ \sum_{e \in \E({\cal T})} \beta_e (x_{\rho(e)}, x_{\tau(e)}) + 
\sum_{e \in \E \backslash \E({\cal T})} \beta_e (x_{\rho(e)}, x_{\tau(e)}) +
\sum_{v \in \V} \alpha_v(x_v).  \]
Adding one edge $e$ of $\G$ in ${\cal T}$ will create a cycle, called a fundamental cycle. In fact, there is a one-to-one correspondence between 
edges $e$ outside the spanning tree and the fundamental cycles ${\cal C}_{e}$ of $\G$ (see \cite{Voloshin}).  So we can rewrite the expression above considering edges corresponding to even cycles and odd cycles, that is, cycles with even or odd number of edges:
\begin{equation} \label{eq: non bipartite}
f(x) \ = \ \sum_{e \in \E({\cal H})} \beta_e (x_{\rho(e)}, x_{\tau(e)}) + 
\sum_{e \in \tilde{\E}} \beta_e (x_{\rho(e)}, x_{\tau(e)}) +
\sum_{v \in \V} \alpha_v(x_v), 
\end{equation}
where ${\cal H}$ is the maximal subgraph of $\G$ containing no odd cycles , ${\cal T} \subseteq {\cal H} \varsubsetneq \G,$ and  
 \[ \tilde{\E} \ = \ \{ \tilde{e} \in \E \backslash \E({\cal T}) \ : \ {\cal C}_{\tilde{e}}  \ {\rm is \  an \ odd \ cycle} \}. \]
${\cal H}$ is bipartite, so it defines a partition $\V = \V_1 \mathbin{\dot{\cup}} \V_2$ so that 
\begin{equation} \label{eq: non bipartite 2}
f(x) \ = \ \sum_{v \in \V_1} \sum_{(u,v) \in  \E({\cal H})} \beta (x_{v}, x_{u}) + 
\sum_{\tilde{e} \in \tilde{\E}} \beta_{\tilde{e}} (x_{\rho(\tilde{e})}, x_{\tau(\tilde{e})}) +
\sum_{v \in \V} \alpha_v(x_v).
\end{equation} 
Let $\tilde{e} = (v_1, v_2) \in \tilde{\E}.$ Notice that $v_1, v_2 \in \V_1$ or $v_1, v_2 \in \V_2,$ otherwise ${\cal C}_{\tilde{e}}$ would be an even cycle. Consider $(u_1, v_1), (u_2, v_2) \in \E({\cal T})$ neighbour edges of $\tilde{e}$.  If $v_1, v_2 \in \V_1,$ the invariance of $\partial f / \partial x_{v_1}$ by permutation of $x_{u_1}$ and $x_{v_2}$ and of $\partial f / \partial x_{v_2}$ by permutation of $x_{u_2}$ and $x_{v_1}$ implies that
\[ \frac{\partial \beta}{\partial x} (x_{v_1}, x_{u_1}) +  \frac{\partial \beta_{\tilde{e}}}{\partial x} (x_{v_1}, x_{v_2}) = \frac{\partial \beta}{\partial x} (x_{v_1}, x_{v_2}) +  \frac{\partial \beta_{\tilde{e}}}{\partial x} (x_{v_1}, x_{u_1}), \]
\[ \frac{\partial \beta}{\partial x} (x_{v_2}, x_{u_2}) +  \frac{\partial \beta_{\tilde{e}}}{\partial y} (x_{v_1}, x_{v_2}) = \frac{\partial \beta}{\partial x} (x_{v_2}, x_{v_1}) +  \frac{\partial \beta_{\tilde{e}}}{\partial y} (x_{u_2}, x_{v_2}). \]
 Hence, for any $a,b$,
 \[\frac{\partial \beta_{\tilde{e}}}{\partial x} (a,b) \ = \   \frac{\partial \beta}{\partial x} (a,b), \ \ \frac{\partial \beta_{\tilde{e}}}{\partial y} (b,a) \ = \   \frac{\partial \beta}{\partial x} (a,b). \] 
 Therefore, $\beta_{\tilde{e}}$ is ${\bf Z}_2$-invariant and $\beta_{\tilde{e}} = \beta.$
If $v_1, v_2 \in \V_2,$ the result follows analogously. Finally, if  $v_1, v_2 \in \V$ and $d(v_1) = d(v_2)$, then at $\bar{x}$ as given in (\ref{eq: x bar}) we have $\partial f / \partial x_{v_1}(\bar{x}) = \partial f / \partial x_{v_2}(\bar{x})$, and using the ${\bf Z}_2$-invariance of $\beta$ we get
\[ d(v_1) \frac{\partial \beta}{\partial x} (a,b) + \alpha^{\prime}_{v_1}(a) \ = \ 
d(v_2) \frac{\partial \beta}{\partial x} (a,b) + \alpha^{\prime}_{v_2}(a), \]
therefore $\alpha_{v_1} = \alpha_{v_2}$ yielding (\ref{eq: eq2 thm admissible functions}).
\hfill \eproof \\

As a consequence of this theorem, a necessary condition for a network function $f : P \to \r$ to be admissible is 
\begin{equation} \label{eq: second derivatives}
\frac{\partial^2 f}{\partial x_{u} \partial x_{v}} \ \equiv \ 0, \ {\rm if}  \ (u,v) \notin \E,
\end{equation}
which captures the `dependency'  relation implied by the network structure (see \cite{De Ville Lerman}). \\

\begin{corollary} \label{cor: admissible functions}
 If $\G$ is a regular graph, then the general form of an admissible function reduces to
\[ f(x) \ = \ \sum_{e \in \E} \phi(x_{\rho(e)}, x_{\tau(e)}), \]
for some smooth ${\bf Z}_2$-invariant function $\phi : \r^{2k} \to \r$.
\end{corollary}
\dem If $\G$ is a regular graph, then self-connection functions are all identical. 

\hfill \eproof

From now on, whenever we label vertices with numbers, their variables shall be indexed with the corresponding numbers. However the subscript of a self-coupling function still denotes the degree of the corresponding vertex.

We illustrate Theorem~\ref{thm: admissible functions} with three graphs: 

\begin{example} \normalfont  \label{ex: admissible functions}
The admissible functions associated with each of the graphs presented in Figure~\ref{fig: three graphs} are given, respectively,  as follows:
\begin{enumerate}
\item $f_{1}(x) = \eta(x_2, x_1) + \eta(x_2, x_3) + \eta(x_2, x_4)$, for some smooth  
$\eta : \r^{2k} \to \r$.
  \item $f_{2}(x) =  \beta(x_1, x_2) + \beta(x_2, x_3) + \beta(x_2, x_4)  + \beta(x_3, x_4) + \alpha_3(x_2) +\alpha_2(x_3) + \alpha_2(x_4)$, for smooth $\alpha_i : \r^k \to \r$, $i=2,3$ and   $\beta: \r^{2k} \to \r$ ${\bf Z}_2$-invariant . 
 \item $f_{3}(x) =  \phi(x_1, x_2) + \phi(x_1, x_5) + \phi(x_1, x_6) + \phi(x_2, x_3) + \phi(x_2, x_4)  + \phi(x_3, x_4) + \phi(x_3, x_5) + \phi(x_4, x_6) + \phi(x_5, x_6)$, with $\phi: \r^{2k} \to \r$ ${\bf Z}_2$-invariant . 
\end{enumerate}
\small{
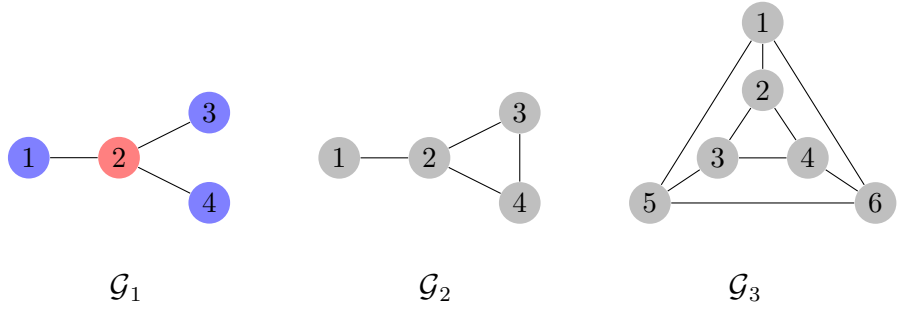
\begin{figure}[H] 
\begin{center}
\tiny{ \begin{tikzpicture}
  [scale=.3,auto=left,every node/.style={circle,fill=black!5}]
  \node[fill=blue!50]  (n1) at (1,9) {\small{1}};
  \node[fill=red!50]  (n2) at (5,9)  {\small{2}};
  \node[fill=blue!50]  (n3) at (9,11)  {\small{3}};
  \node[fill=blue!50]  (n4) at (9,7) {\small{4}};
  
  \foreach \from/\to in {n1/n2,n2/n3,n2/n4}
    \draw (\from) -- (\to);
\end{tikzpicture}  \hspace*{1cm}
\begin{tikzpicture}
  [scale=.3,auto=left,every node/.style={circle,fill=black!25}]
  \node (n1) at (1,5) {\small{1}};
  \node (n2) at (5,5)  {\small{2}};
  \node (n3) at (9,7)  {\small{3}};
  \node (n4) at (9,3) {\small{4}};
  
  \foreach \from/\to in {n1/n2,n2/n3,n2/n4,n3/n4}
    \draw (\from) -- (\to);
\end{tikzpicture} \hspace*{1cm}
\begin{tikzpicture}
  [scale=.3,auto=left,every node/.style={circle,fill=black!25}]
  \node (n1) at (9,4) {\small{1}};
  \node  (n2) at (9,1)  {\small{2}};
  \node  (n3) at (7,-2)  {\small{3}};
  \node  (n4) at (11,-2) {\small{4}};
 \node  (n5) at (4,-4)  {\small{5}};
  \node (n6) at (14,-4) {\small{6}};
  
  \foreach \from/\to in {n1/n2,n2/n3,n2/n4,n3/n4,n3/n5,n4/n6,n5/n6, n1/n5, n1/n6}
    \draw (\from) -- (\to); 
\end{tikzpicture}}
\end{center}
\begin{center} $\hspace*{3.7cm}\G_1$ \hfill $\G_2$ \hfill $\G_3$ \hspace*{4cm} \end{center}
\caption{\small{Graph $\G_i$ corresponding to the admissible function $f_i$ given in Example~\ref{ex: admissible functions}, $i=1,2,3$.}} \label{fig: three graphs}
\end{figure}}
We just notice that we first use Theorem~\ref{thm: admissible functions} to write  $f_1$ as
\[ f_1(x)  = \beta(x_2, x_1) +  \beta(x_2, x_3) + \beta(x_2, x_4) + 
\alpha_1(x_2) + \alpha_1(x_1) + \alpha_1(x_3) + \alpha_1(x_4),\]
and then define $\eta(a,b) = \beta(a,b) + \alpha_1(b).$ 
\end{example}
We give below another example of a bipartite graph, with a larger number of cells:

\begin{example} \label{ex: another bipartite}
{\rm We consider a bipartite graph for which $d(u) \neq d(v)$ if $u$ and $v$ have distinct colours, see  Figure~\ref{fig: another bipartite graph}. So an admissible function for this graph is of the form
(\ref{eq: eq1 thm admissible functions}) and it is given by
\begin{align*}
 f(x) =  &   \sum_{i=1,3,5,7,8,9} \beta(x_{2}, x_{i}) + \sum_{j=3,5,7,8,9} \beta(x_{4}, x_{j})  +  \sum_{k=3,8,9,10} \beta(x_{6}, x_{k}) + \alpha_1(x_{1}) +  \alpha_1(x_{{10}}) + \\
& +  \alpha_2(x_{5})+ \alpha_2(x_{7}) + \alpha_3(x_{3}) + \alpha_3(x_{8}) + \alpha_3(x_{9}) + \gamma_4 (x_{6})   +  \gamma_5(x_{4}) + \gamma_6(x_{2}).
\end{align*}

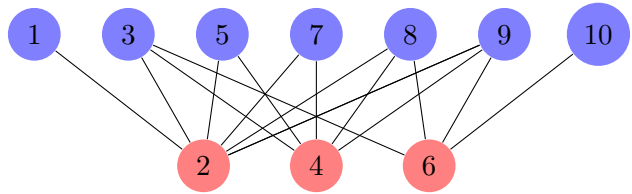
\begin{figure}[H] 
\begin{center}
\begin{tikzpicture}
  [scale=.25,auto=left,every node/.style={circle,fill=black!5}]
  \node[fill=blue!50] (n1) at (-2,9) {\small{1}};
  \node[fill=blue!50] (n3) at (3,9)  {\small{3}};
  \node[fill=blue!50] (n5) at (8,9)  {\small{5}};
  \node[fill=blue!50] (n7) at (13,9) {\small{7}};
   \node[fill=blue!50] (n8) at (18,9) {\small{8}};
  \node[fill=blue!50] (n9) at (23,9)  {\small{9}};
  \node[fill=blue!50] (n10) at (28,9)  {\small{10}};
  \node[fill=red!50] (n2) at (7,2) {\small{2}};
  \node[fill=red!50] (n4) at (13,2)  {\small{4}};
  \node[fill=red!50] (n6) at (19,2)  {\small{6}};

  \foreach \from/\to in {n1/n2,n2/n3,n3/n4, n3/n6, n5/n2, n5/n4, n7/n2, n7/n4, n8/n2, n8/n4, n8/n6, n9/n2, n9/n2, n9/n6, n10/n6, n9/n4}
    \draw (\from) -- (\to);
 \end{tikzpicture}
\end{center}
\caption{\small{A bipartite graph with 10 cells.}} \label{fig: another bipartite graph}
\end{figure}
}
\end{example}


\section{Critical points on regular graphs}  \label{sec:reg graphs}

In this section we direct attention to regular graphs, so admissible functions are of the form
\begin{equation} \label{eq: fn regular graph}
f(x) =  \sum_{ e \in \E} \phi(x_{\rho(e)}, x_{\tau(e)}), 
\end{equation}
with $\phi: \r^{2k} \to \r$ ${\bf Z}_2$-invariant. Our aim is to discuss when the nature of certain types of critical points of $f$ can be deduced from the associated critical points of the coupling function $\phi$.  Totally synchronous critical points and  2-colour patterns of critical points are considered, and these are given in  Subsection~\ref{subsec: tot synchr regular} and in Subsection~\ref{subsec: 2-colour regular}, respectively.  Minima of one function are detected from minima of the other depending on the graph architecture. It is clear that our analysis is distinct for these two types of critical points  in the class of admissible functions. Critical configurations with three or more colours must be discussed under a different approach,  using results in  \cite{ADM}.

\subsection{Totally synchronous patterns}  \label{subsec: tot synchr regular} 

For an admissible function of the form (\ref{eq: fn regular graph}), a point $x^0 = (x_0, \ldots, x_0)$ is a critical point of $f$ if, and only if, $(x_0, x_0)$ is a critical point of its coupling $\phi$, even when $\phi$ has no permutation invariance. In this subsection we assume that the graph is regular and that the configuration space of each cell is one-dimensional ($k=1$). In this case,  the coupling is ${\bf Z}_2$-invariant, and so  the nature of such critical points given by the Hessian criterion is then a 2-parameter problem, the parameters being given by the second partial derivatives  of  $\phi$  at $(x_0, x_0) \in \r^2$. \\

We shall use below subscripts 1 and 2 for the partial derivatives of $\phi$ with respect to the first and second variable respectively.  In what follows we denote $\alpha = \phi_{11}(x_0, x_0)$ and $\beta = \phi_{12}(x_0, x_0)$. The Hessian of $f$
at $x^0$ is the matrix $H = (\alpha + \beta) d {\rm I} - \beta {\cal L}$, where ${\rm I}$ is the identity matrix and ${\cal L} = D + A$ is the graph Laplacian, where $D$ is the degree matrix and $A$ is the adjacency matrix of $\G$ (\cite{Voloshin}). Notice that the eigenvalues of $H$ are
\begin{equation} \label{eq:evalues of H}
 \mu_i = (\alpha + \beta) d - \beta \lambda_i, 
 \end{equation}
where $0 = \lambda_1 < \lambda_2 \leq \ldots \leq \lambda_n = \lambda$ are the eigenvalues of ${\cal L}$. 

\begin{proposition} \label{prop: tot synch reg graphs}
Let $\G$ be a $d$-regular graph and consider an admissible function on $\G$ of the form
(\ref{eq: fn regular graph}). Then generically a totally synchronous critical point 
$x^0 = (x_0, \ldots, x_0)$ is a minimum of $f$ if, and only if, 
\[ \phi_{11}(x_0, x_0) + \phi_{12}(x_0, x_0)  > 0, \ \ 
 \phi_{11}(x_0, x_0) + (1 - \frac{\lambda}{d})\phi_{12}(x_0, x_0) >0, \] 
where $\lambda$ is the greatest eigenvalue of the Laplacian of $\G$. 
\end{proposition}
\dem
The Hessian  $H$ is generically nondegenerate, in which case $x^0$ is a minimum if, and only if, 
\[ \alpha > \frac{-d + \lambda_i}{d} \beta, \]
for all  $i = 1, \ldots, n.$
It is well-known that $\lambda \leq 2d,$ where equality holds if, and only if, $\G$ is bipartite. Hence,
for any $i=1, \ldots, n$,  we have $(-d + \lambda_i)/d \in [-1,1]$, and so a minimum occurs if, and only if, 
\begin{equation} \label{eq: inequalities}
\alpha + \beta > 0 \  \ {\rm and}  \  \ \alpha + (\frac{d-\lambda}{d}) \beta > 0. 
\end{equation}
\hfill \eproof

Inequalities (\ref{eq: inequalities}) are precisely the Hessian criterion  for a minimum of  $\phi$  in (\ref{eq: fn regular graph}) if $\lambda = 2d$. So, for the generic case when this criterion applies, we have:

\begin{corollary} \label{cor: tot synch reg graphs1}
An  admissible  function $f$ defined on a regular graph, given by (\ref{eq: fn regular graph}),  has a minimum at $x^0 = (x_0, \ldots, x_0)$ if
 the coupling function $\phi$ has a minimum at $(x_0, x_0)$. The converse holds if, and only if, $\G$ is bipartite.
\end{corollary}
  
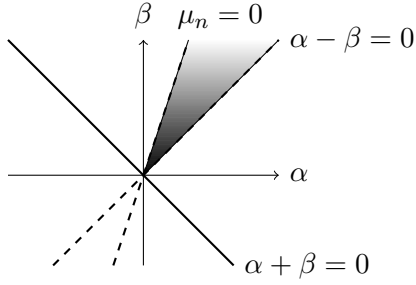
\begin{figure} 
\begin{center}
\begin{tikzpicture}[scale=0.6]
\draw[top color=white,bottom color=black] (3,3) -- (0,0) -- (1,3) ; 

\foreach \coordinate/\label/\pos in {{(2,-2)/{\small $\alpha+\beta=0$}/right},{(3,3)/{\small $\alpha-\beta=0$}/right},{(1.7,3)/{\small $\mu_n=0$}/above}}
\node[\pos] at \coordinate {\label};
\draw[->] (-3.0,0) -- (3.0,0) node[right] {{\small $\alpha$}} coordinate(x axis); 
\draw[->] (0,-2) -- (0,3.0) node[above] {{\small $\beta$}} coordinate(y axis); 
\draw[black, thick, dashed] (-2,-2) -- (3.0,3.0);
\draw[black, thick, dashed] (-2/3,-2) -- (1,3);
\draw[black, thick] (-3.0,3.0) -- (2,-2);
\end{tikzpicture} 
\end{center}
\caption{{\small The wedge region where the converse of Corollary~\ref{cor: tot synch reg graphs1} fails.}} \label{fig: wedge}
\end{figure}  
Looking at the plane of the parameters $\phi_{11}(x_0, x_0),$ $\phi_{12}(x_0, x_0)$,   the converse in the corollary above fails precisely in the wedge
\begin{equation} \label{eq: wedge}
 \phi_{11}(x_0, x_0) - \phi_{12}(x_0, x_0)  <  0, \  \  
\phi_{11}(x_0, x_0) + (1 - \frac{\lambda}{d})\phi_{12}(x_0, x_0) >0, 
\end{equation}
see Figure~\ref{fig: wedge}. Now, for a $d$-regular graph, we have the lower bound for the greatest eigenvalue 
$\lambda,$  
\[ \lambda \geq d + 1, \]
where equality holds if, and only if, the graph is complete \cite{Zhang}. This means that we can measure the set in parameter space where minima of $f$ are not in one-to-one correspondence with minima of the coupling function $\phi$. Its `size'  is determined by the topology of the graph:

\begin{corollary} \label{cor: tot synch reg graphs2}
For regular graphs, on the plane of the two parameters $\phi_{11}(x_0,x_0)$ and
$\phi_{12}(x_0,x_0)$ the wedge determined by (\ref{eq: wedge}) vanishes if, and only if, the graph is bipartite. It is as large as possible if, and only if, the graph is complete.
\end{corollary} 

More generally, for any regular graph we can use  (\ref{eq:evalues of H})  to track down the nature of a fully synchronous critical point $x^0$ of $f$ for different values of  $\alpha$ and $ \beta$. Equivalently, we detect the degree of instability of  $x^0$ thought of as an equilibrium of the  vector field  $- \nabla f.$  In fact,  in the open region $(\alpha+\beta)(\alpha - \beta) >0$, $x^0$ is a minimum of $f$  if $\alpha>0$ and it is maximum if $\alpha < 0$.  If $n$ is the number of cells in $\G$, we have that in the region $(\alpha+\beta)(\alpha - \beta) < 0$ are the lines on which  $\mu_j=0,$ for $1 < j  < n$, whose equations are
\[ \alpha \ = \ \frac{-d+\lambda_j}{d} \beta, \]
with algebraic multiplicity of $\mu_j$ given by the algebraic multiplicity of $\lambda_j$ for $L$. On this line and for $\beta >0,  \mu_k > 0$ if $1 \leq k < j$ and  $ \mu_k < 0$ if $j < k < n$;  for 
$\beta <0$, $\mu_k < 0$ if $1 \leq k < j$ and $\mu_k >0$ if $j < k < n$. Also,  the line
$\alpha + \beta =0$ corresponds to the case for which $H$ is a $\beta$-weighted Laplacian, that is, a weighted Laplacian (see \cite{Voloshin}) with all weights equal to $\beta$. So 
$\mu_1=0$ and a minimum of $f$ can only occur for negative values of $\beta$.  Finally, on the line 
$\alpha- \beta = 0$, for positive (negative) values of $\alpha$ we have $\mu_j > 0 (<0)$, for $ 1 \leq j \leq n$ unless the graph is bipartite, in whose case all eigenvalues are still positive (negative) except $\mu_n$, which vanishes. \\

With the two corollaries above in mind, it becomes natural to look for a class of nonregular graphs for which
the wedge vanishes, that is, for which the converse of Corolllary~\ref{cor: tot synch reg graphs1}
does hold. As we see in the next proposition, this is the case for the complete bipartite graph $K_{m,n}$, i.e., the bipartite graph where every vertex of one set of vertices $\V_1$ is attached to every vertex of the other set $\V_2, $ with $|\V_1| = m, |\V_2| = n$.  If $m=n$, $K_{m,n}$ is regular and it fits in Corollary~{\ref{cor: tot synch reg graphs1}, so in the next result we are interested in the case  $ m \neq n$.

\begin{proposition} \label{prop: tot synch reg graphs3}
For distinct $m,n \geq 2,$  an  admissible  function $f$ given in (\ref{eq: fn regular graph}) defined on the complete bipartite graph $K_{m,n}$ has a minimum at $x^0 = (x_0, \ldots, x_0)$ if, and only if, 
 the coupling $\phi$ has a minimum at $(x_0, x_0)$.
\end{proposition}

\dem
The Hessian matrix of $f$ at $x^0$ is given by
\[ H(x^0) \  = \   \left( \begin{array}{c|c} 
n \alpha {\rm I}_{m} & \beta  B^t   \\ \cline{1-2} 
\beta B & m \alpha {\rm I}_n \end{array} \right), \]
where $\alpha = \phi_{11}(x_0, x_0),$  $\beta = \phi_{12}(x_0, x_0)$ and $B$ is the
$n \times m$-matrix in the adjacency matrix 
\begin{equation} \label{eq: B-adjacency}
A \ = \ \left( \begin{array}{c|c}
0 & B^t  \\ \cline{1-2}
B& 0 
\end{array} \right)
\end{equation} 
of $K_{m,n}$, with all entries equal to 1.  We shall now describe  the  elements in the spectrum
Spec$(H(x^0))$. 
 For $\lambda \neq n \alpha$, we have that $\det(H(x^0) - \lambda {\rm I}_{m+n}) = 0$ if, and only if,
\[ \det\Bigl((m \alpha - \lambda)( n \alpha - \lambda) {\rm I}_n - \beta^2B B^{t} \Bigr)  \ = \ 0. \]
The matrix in the expression above can be rewritten as
\[ {\rm I}_n \lambda^2- (m+n) \alpha  {\rm I}_n \lambda + mn \alpha^2 {\rm I}_n - \beta^2B B^t  \ = \ 
(M_{+} - \lambda {\rm I}_n) (M_{-} - \lambda {\rm I}_n), \]
where
\[ M_{\pm} \ = \   \ \frac{1}{2} \bigl[    (m + n) \alpha    {\rm I}_n \pm 
{\Bigl( (m-n)^2 \alpha^2 {\rm I}_{n} + 4 \beta^2 B B^t  \Bigr)}^{1/2} \bigr],  \]
which are matrices with real entries since $B B^{t}$ is a positive semidefinite matrix. We now investigate the eigenvalues of these two matrices distinct from $n \alpha$. We notice that the eigenvalues of $B B^t$ are given by $\mu = \zeta^2$, for each eigenvalue  $\zeta$ of the adjacency matrix $A$. But it is well-known that $\zeta = 0, - \sqrt{mn}, \sqrt{mn}$. Hence, the eigenvalues of $M_{+}$ and $M_{-}$ are 
\[ \frac{1}{2} \bigl[ (m+n) \alpha \pm \Bigl( (m-n)^2 \alpha^2 + 4 \mu \beta^2 \Bigr)^{1/2} \bigr] ,\] 
for each of the two values of $\mu$. Therefore, the eigenvalues distinct from $n \alpha$ are  $m \alpha$ and  
\[ \frac{1}{2} \bigl[ (m+n) \alpha \pm \Bigl( (m-n)^2 \alpha^2 + 4 mn  \beta^2 \Bigr)^{1/2} \bigr]. \] 
For $\lambda \neq m \alpha$ we proceed analogously as above for $\lambda \neq n \alpha$ to conclude that  $n \alpha$ is also an eigenvalue of $H(x^0)$. And now it is straightforward to check that for any $m$ and $n$  these eigenvalues are all positive if, and only if,   $(\alpha + \beta),  (\alpha - \beta)> 0$.  
\hfill \eproof

\subsection{2-colour patterns }  \label{subsec: 2-colour regular} 

In this subsection we investigate critical points $\bar{x}$ of admissible functions (\ref{eq: fn regular graph}) on regular cell graphs $\G$ which are not totally synchronous and which come from critical points $(x_{0}, y_{0})$ of the coupling function with $x_0 \neq y_0$.  So,  for $\V = \{1,2, \ldots, n\},$ we are considering
\begin{equation} \label{eq: x bar 2} 
\bar{x} = (x_1, \ldots, x_n), \ \   x_i = x_0 \ {\rm or} \  y_0, \ \ i=1, \ldots, n. 
\end{equation}
We start by showing that no such critical points are expected  unless the graph is bipartite, even if it is non-regular:

\begin{proposition} \label{prop: generic 2-colour}
Let $\G = (\V, \E)$ be a cell graph. For a generic coupling function $\phi$ a point $\bar{x} \in \r^{kn}$ as in (\ref{eq: x bar 2}) can be a critical point of an admissible function of the form
\begin{equation} \label{eq: fn non regular graph}
f(x) \ = \ \sum_{e \in \E}  \phi(x_{\rho(e)}, x_{\tau(e)}), 
\end{equation}
$\phi : \r^{2k} \to \r,$  only  if $\G$ is bipartite. In this case, for $\V =  \V_1 \mathbin{\dot{\cup}} \V_2$, $ \bar{x}$ is given by 
\begin{equation} \label{eq: critical point}
x_{i} = \left\{ \begin{array}{ll} 
x_0,  &  i \in \V_1  \\
y_0, & i \in \V_2. \end{array} \right. 
\end{equation}
\end{proposition}
\dem We may assume that the coupling function $\phi$ is ${\bf Z}_2$-invariant; otherwise, the graph is already bipartite, by Theorem~\ref{thm: admissible functions}. Suppose that $(x_0, y_0)$ is a critical point of $\phi$ with the generic condition that neither are $(x_0, x_0)$ nor  $(y_0, y_0).$ Fix $j$, $1 \leq j \leq n$. If $x_{j} = x_0,$ let $p$ be the number of vertices $\ell \in I(j)$ such that $x_{\ell} = x_0$ and  $q$ be the number of vertices $\ell \in I(j)$ such that $x_{\ell} = y_0$. Then for any $i$, $1 \leq i \leq k$,  
\[ \frac{\partial f}{\partial x_j^i} (\bar{x}) \ = \sum_{{\ell} \in I(j)} \frac{\partial \phi}{\partial x_j^i} (x_j, x_\ell) =  \ p \ \frac{\partial \phi}{\partial x_j^i} (x_0, x_0) + q \ \frac{\partial \phi}{\partial x_j^i}(x_0, y_0) =  p \ \frac{\partial \phi}{\partial x_j^i} (x_0, x_0). \]
By the permutation invariance of $\phi$,  $\partial \phi/\partial x_j^i (x_0, x_0) \neq 0$
for some $i$, $1 \leq i \leq k.$ Hence, $p=0$, and therefore $x_\ell = y_0$ for all $l \in I(j)$.   Analogously, if $x_j = y_0$, then $x_\ell = x_0$ for all $\ell \in I(j)$. This defines a partition of $\V$ so that $\G$ is bipartite and the components of $\bar{x}$  must be given by ({\ref{eq: critical point}).
\hfill \eproof

As a consequence of the proposition above,  the assignment of  one distinct colour for each $x_0$ and $y_0$ provides a standard 2-colouring on the graph, that is, no two vertices sharing the same edge have the same colour.  Let us also notice that in the case when the coupling function is invariant by permutation, then 2-colour critical points of $f$ come in pairs. \\

If $\G$ is a regular bipartite graph, then it has an even number of vertices,  $n=2m$, $m \geq 1$, with  
$|\V_1| = |\V_2| = m$. In fact, just impose the regularity condition on the  well-known equality for
 bipartite graphs,  
\[\sum_{u \in \V_1} \deg(u) = \sum_{v \in \V_2} \deg(v).\]
So here we shall investigate  2-colour patterns constrained to the set of $(d,m)$-graphs, as we define next:

\begin{definition} \label{def: d-m graph}
A connected graph $\G = (\V, \E)$ is a $(d,m)$-graph if it is a $d$-regular bipartite graph
with $| \V_1| = |\V_2| = m$.
\end{definition}

We refer to \cite{Meringer1, Meringer2} for the number of existing  $(d,m)$-graphs when $d=3,4,5$ and $3 \leq m \leq 16.$ 
Figure~\ref{fig:  cubic graph} illustrates the cubic graph, which is a $(3,4)$-graph.

\begin{figure}[H] 
\tiny{\begin{center}
\begin{tikzpicture}
 [scale=.15,auto=left,every node/.style={circle,fill=black!15}]
 \node[fill=blue!50] (n4) at (1,1) {\tiny{4}};
 \node[fill=red!50] (n7) at (1,13)  {\tiny{7}};
 \node[fill=red!50] (n6) at (5,5)  {\tiny{6}};
 \node[fill=blue!50] (n1) at (5,9) {\tiny{1}};
 \node[fill=blue!50] (n2) at (9,5) {\tiny{2}};
 \node[fill=red!50] (n5) at (9,9)  {\tiny{5}};
 \node[fill=blue!50] (n3) at (13,13)  {\tiny{3}};
\node[fill=red!50] (n8) at (13,1)  {\tiny{8}};
 \foreach \from/\to in {n4/n8,n4/n7,n4/n6, n2/n8, n2/n5, n2/n6,n1/n6,n1/n5,n1/n7,n3/n7,n3/n5,n3/n8}  \draw (\from) -- (\to);
\end{tikzpicture}
\end{center}}
\caption{\small{The cubic graph is a $(3,4)$-graph.}} \label{fig: cubic graph}
\end{figure}
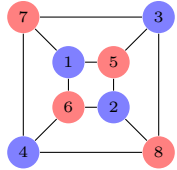

In the next two propositions we assume the variable on each cell is one-dimensional.

\begin{proposition} \label{prop: invariant spectrum}
The spectrum of a 2-colour critical point of an admissible function depends only on the coupling function and the integers $d$ and $m$.
\end{proposition}

\dem
Given any $(d,m)$-graph $\G =(\V, \E) $, $\V =  \V_1 \mathbin{\dot{\cup}} \V_2$, use 
$\{1,2, \ldots, m\}$ to label elements in $\V_1$ and $\{m+1,m+2, \ldots, 2m\}$ to label elements in 
$\V_2.$ For any admissible function (\ref{eq: fn regular graph}) on $\G$, its Hessian at the 2-colour
critical point  $\bar{x}$ given by (\ref{eq: critical point}) is of the form
\begin{equation} \label{eq: d-m Hessian}
H(\bar{x}) \  = \  \left( \begin{array}{c|c}
d \alpha {\rm I}_m & 0 \\ \cline{1-2}
0 & d \gamma {\rm I}_m 
\end{array} \right)  + \beta A, 
\end{equation} 
where $\alpha = \phi_{11}(x_0, y_0)$,  $\beta = \phi_{12}(x_0,  y_0)$, $\gamma = \phi_{11}(y_0, x_0)$ 	
and $A$ is the adjacency matrix of $\G$, which is of the form
\begin{equation} \label{eq: C-adjacency}
A \ = \ \left( \begin{array}{c|c}
0 & C^t  \\ \cline{1-2}
C& 0 
\end{array} \right),  
\end{equation} 
for an order-$m$ matrix $C.$ 

Any other $(d,m)$-graph $\tilde{\G}$ is given by permuting edges of $\G$, so there exists an 
order-$m$ permutation matrix $P$ such that its adjacency matrix is given by

\[ \tilde{A} \ = \ \left( \begin{array}{c|c}
0 & (PC)^t  \\ \cline{1-2}
PC & 0 
\end{array} \right).  \] 
But the Hessian matrix for $\tilde{\G}$ at $\bar{x}$ is
\[ \tilde{H}(\bar{x}) \ = \  \left( \begin{array}{c|c}
d \alpha {\rm I}_m & 0 \\ \cline{1-2}
0 & d \gamma {\rm I}_m 
\end{array} \right)  + \beta \tilde{A}, \]
which is similar to (\ref{eq: d-m Hessian}): just take 
\[ \left( \begin{array}{c|c}
P & 0 \\ \cline{1-2}
0 & {\rm I}_m 
\end{array} \right) \]
as a similarity matrix.
\hfill \eproof

The next result  is the analogous of Proposition~\ref{prop: tot synch reg graphs3} for $(d,m)$-graphs at 2-colour pattern configurations. 

\begin{proposition} \label{prop: 2-colour reg bipartite}
If $\G$ is a $(d,m)$-graph,  then a 2-colour pattern $\bar{x}$  in $\G$ as in (\ref{eq: critical point})  is a minimum configuration of  an  admissible  function given in (\ref{eq: fn regular graph}) defined on $\G$ if, and only if,   $(x_0, y_0)$ is a minimum of 
 its coupling function.
\end{proposition}

\dem
The Hessian matrix of $f$ at ${\bar{x}}$ is given by (\ref{eq: d-m Hessian})-(\ref{eq: C-adjacency}), and
its characteristic polynomial  $p_{H(\bar{x})}$ is given by the product of two characteristic polynomials
 \[ p_{H(\bar{x})} \ = \  p_{N_{+}}  \  p_{N_{-}} ,   \]
where
\[ N_{\pm} \ =  \ \frac{1}{2} \bigl[   d (\alpha + \gamma)    {\rm I}_m \pm 
{\Bigl(d^2 (\alpha-\gamma)^2 {\rm I}_{m} + 4 \beta^2 C C^t  \bigr)}^{1/2} \bigr],  \]
which are matrices with real entries since $C C^{t}$ is a positive semidefinite matrix. Hence,
\[ {\rm Spec}(H(\bar{x})) = {\rm Spec}(N_{+}) \cup {\rm Spec}(N_{-}). \]
Now the least element of  ${\rm Spec}(H(\bar{x}))$ is the minimum of ${\rm Spec}(N_{-})$, which is given by
\[ \xi_{N_{-}} = \frac{1}{2}  \Bigl(d(\alpha +\gamma) - \Bigl( d^2(\alpha - \gamma)^2
+ 4 \nu \beta^2 \Bigr)^{1/2} \Bigr), \]
where $\nu$ is the greatest eigenvalue of $C C^{t}$. But $\nu = \lambda^2$, where $\lambda$ is the largest eigenvalue of the adjacency matrix of $\G$. Since $\G$ is $d$-regular and bipartite, it follows that $\lambda = d$. Therefore,
\[ \xi_{N_{-}} = \frac{d}{2}  \Bigl((\alpha +\gamma) - \Bigl( (\alpha - \gamma)^2
+ 4  \beta^2 \Bigr)^{1/2} \Bigr), \]
which is positive if, and only if, $\alpha > 0$ and  $\alpha \gamma - \beta^2 >0$. 
\hfill \eproof

We end this section with a  remark regarding graphs in general, not necessarily regular.

\begin{remark} \normalfont \label{rmk:coexistence}
As a consequence of Proposition~\ref{prop: generic 2-colour}, unless there are additional symmetries other than possible permutation invariance of the coupling function, occurrence of both 2-colour patterns and totally synchronous critical points
$(z_0, \ldots, z_0) $ is expected for admissible functions (\ref{eq: fn non regular graph}) only if $z_0 \neq x_0, y_0$.  However, under extra symmetry constraints, these two types of configurations can co-exist for  $z_0 = x_0$ or $y_0$ . In fact, these can be critical configurations, among others, for example when the admissible function has an invariance under  the circle group  ${\bf S}^1$. In this case these two can even both be minima of $f$. We shall return to this point  in the next section. A  study of a particular ${\bf S}^1$-invariant function on  a ring is carried out in Subsection~\ref{subseq: minimum on a ring}.
\end{remark}


\section{${\bf S}^1$-invariant admissible functions}  \label{sec:S1-invariant admissible functions}

In this section we analyse  possible critical configurations of a class of admissible functions for which the interest lies in the behaviour of phase differences between  planar unit vector states  of coupled cells. For each cell $i$ in the graph we can identify the vector state with its angle $\theta_i$ with respect to a fixed direction, which we take to be the vertical direction. In this setting,  for a network represented by a graph $\G$ of $n$ cells,  the  configuration manifold is the $n$-torus $\T^n = \S^1 \times \cdots \times \S^1.$  
In this case the admissible functions on $\cN = (\G, \T^n)$ 
are additionally invariant under the translational diagonal action of $\S^1,$
\[ h(\theta_1 + \psi,\ldots  \theta_n + \psi) = h(\theta_1, \ldots, \theta_n), \ \forall \psi \in \S^1. \]
Imposing this condition  on the general  form  given by Theorem~\ref{thm: admissible functions} yields
\begin{equation} \label{eq: self-connection}
 \alpha_{d(v)} \equiv 0,  \  \forall v \in \V, 
 \end{equation}
so an $\S^1$-invariant admissible function has only trivial 
self-connections. In addition,  the coupling  $\beta$ in (\ref{eq: eq1 thm admissible functions}) or (\ref{eq: eq2 thm admissible functions}) is of the form
\begin{equation} \label{eq: S1-invariant beta}
 \beta(x , y) = \delta (x-y), 
\end{equation}
for some smooth function $\delta : \S^1 \to \r$. If $\delta$ admits an extension to $\r$ then this must be an even function for the case  (\ref{eq: eq2 thm admissible functions}).  
We remark  that (\ref{eq: self-connection}) is a consequence only of the $\S^1$-invariance, so it also holds for ${\bf S}^1$-invariant functions defined on non-regular cell graphs (see Corollary~\ref{cor: admissible functions}). \\

Functions of this type appear in many applications, including the Kuramoto model (see \cite{Bronski et al, MS}) and the Antiferromagnetic XY model (see \cite{Korshunov, Lee, Walter Chatelain}). We briefly describe the second, which models alignments of spins in magnetic materials. The spins are distributed in a planar lattice, the `site' of each being each vertex of the lattice. The source of the general behaviour of magnetic materials is the spontaneous, parallel alignment of neighbouring spins. This type of magnetism is called {\it ferromagnetism}. An equally common yet very different form of magnetism is {\it antiferromagnetism}. Like ferromagnetism, this arises due to the spontaneous alignment of neighbouring spins, but the alignment of spins in antiferromagnetic materials is anti-parallel (in an up-down configuration). 
The total energy of this spin system, the Hamiltonian, assumes  a particular form of an admissible function defined on the lattice, with coupling functions typically being given by the cosine function and the sum related to couplings is over nearest neighbour spin sites.

\subsection{The Hessian}

Let us first consider the case of nearest-neighbour coupling  of $n$ identical cells in a ring, so
the  cell graph $\G$ is 
a regular $n$-sided polygon.  The admissible function is then given by
\begin{equation} \label{eq: function on a ring}
h(\theta) \ = \ \sum_{i=1}^{n} \delta(\theta_i -  \theta_{i-1}),
\end{equation}
where $\theta_{0} = \theta_n$ and $\delta$ is  even. 

Critical points of $h$ are points $\theta = (\theta_1, \ldots, \theta_n)$ such that
\begin{equation} \label{eq: cr pts on a ring}
\delta^{\prime} (\theta_{i+1} - \theta_{i} )= \delta^{\prime}(\theta_{i} - \theta_{i-1}), \  1 \leq i \leq n-1.
\end{equation}
The Hessian matrix $H$ of $h$ has its entries given as follows:   for $1 \leq i \leq n,$
\[ \frac{\partial^2 h}{\partial \theta_i^2} \ = \ \delta^{\prime \prime}(\theta_{i+1} - \theta_{i}) + \delta^{\prime \prime}(\theta_i - \theta_{i-1}),  \]

\[  \frac{\partial^2 h}{\partial \theta_i \partial \theta_j}  \ = \ \left\{ \begin{array}{ll}
-\delta^{\prime \prime}(\theta_i - \theta_j), &  \ j = i+1, \ i-1 \\
0, & \ {\rm otherwise}. 
\end{array} \right. \]

From this particular case of a ring  it is straightforward to see how the expression of the Hessian 
generalizes for any cell graph $\G$ of $n$ vertices: If $\G$ is bipartite (and $\delta$ may not be even), from the partition $\V = \V_1 \mathbin{\dot{\cup}} \V_2$ of the set of vertices the admissible function is of the form 
\begin{equation} \label{eq: bipartite S1-inv}
h(\theta_1, \ldots, \theta_n) \ =  \sum_{i \in \V_k, \ j \in I(i)} \delta(\theta_i - \theta_j), 
\end{equation}
for arbitrary $k =1$ or 2. Set $l = 1$ or 2, $l \neq k.$ Then the critical points $\theta$ of $h$ have coordinates that solve the equations 
\begin{align} \label{eq: S1 critical bipartite}
\nonumber \\
 \sum_{j \in I(i)} \delta^{\prime}(\theta_i - \theta_{j}) \ = 0, \ \ \ {\rm for} \ \   i \in \V_k \\ \nonumber
 \sum_{j \in I(i)} \delta^{\prime}(\theta_j - \theta_{i}) \ = 0, \ \ \ {\rm for} \ \   i \in \V_l. 
 \end{align}
The Hessian in this case is as follows: for  $1 \leq i \leq n,$
\[ \frac{\partial^2 h}{\partial \theta_i^2}  \ = \left\{ \begin{array}{ll} 
\sum_{j \in I(i)} \delta^{\prime \prime}(\theta_i - \theta_{j}), & \ \ {\rm if} \ \   i \in \V_k \\
\sum_{j \in I(i)} \delta^{\prime \prime}(\theta_{j} - \theta_{i}), & \ \ {\rm if}  \ \ i \in \V_l
\end{array} \right. \]

\[    \frac{\partial^2 h}{\partial \theta_i \partial \theta_j} \ = \ \left\{ \begin{array}{ll}
-\delta^{\prime \prime}(\theta_i - \theta_j), &  \ \ {\rm if} \ \ i \in \V_k, \ \ j \in I(i), \ j \neq i \\
-\delta^{\prime \prime}(\theta_j- \theta_i), &  \ \ {\rm if} \ \ i \in \V_l, \ \ j \in I(i), \ j \neq i \\
0, & \ \ {\rm otherwise}. 
\end{array} \right. \]
If $\G$ is not bipartite,  then the function is of the form
\begin{equation} \label{eq: S1-inv function}
 h(\theta_1, \ldots, \theta_n)  \  = \   \sum_{(i,j) \in \E} \delta(\theta_i - \theta_j), 
 \end{equation}
where $\delta$ is even and the equations for critical points as well as the Hessian entries above simplify in the obvious way.  

\begin{example}\normalfont  \label{ex: bipartite S1 inv}
Consider the graph in Figure~\ref{fig: bipartite graph}. An admissible function defined on this graph is of the form (\ref{eq: bipartite S1-inv}) and the Hessian matrix is given by  
{\small \[ H \ = \ \left( \begin{array}{ccccccc}
a_{21} & - a_{21} & 0 & 0 & 0 & 0 & 0  \\
- a_{21}  & a_{21} + a_{23} & -a_{23} & 0 & 0 & 0 & 0 \\
0 & -a_{23} & a_{23} + a_{43} + a_{53} & - a_{43}  & -a_{53} & 0 & 0 \\
0 & 0 & -a_{43} & a_{43} + a_{46} & 0 & - a_{46} & 0 \\ 
0 & 0 & -a_{53} & 0 & a_{53} + a_{57} & 0 & -a_{57} \\
0 & 0 & 0 & -a_{46} & 0 & a_{46} & 0 \\
0 & 0 & 0 & 0 & -a_{57} & 0 & a_{57} 
\end{array} \right), \] }
where   $a_{ij}= \delta^{\prime \prime}(\theta_i - \theta_j)$. 

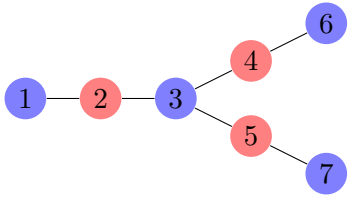
\begin{figure}[H] 
\begin{center}
\tiny{ \begin{tikzpicture}
  [scale=.25,auto=left,every node/.style={circle,fill=black!15}]
  \node[fill=blue!50] (n1) at (1,9) {\small{1}};
  \node[fill=red!50] (n2) at (5,9)  {\small{2}};
  \node[fill=blue!50] (n3) at (9,9)  {\small{3}};
  \node[fill=red!50] (n4) at (13,11) {\small{4}};
  \node[fill=red!50] (n5) at (13,7) {\small{5}};
  \node[fill=blue!50] (n6) at (17,13)  {\small{6}};
  \node[fill=blue!50] (n7) at (17,5)  {\small{7}};
\foreach \from/\to in {n1/n2,n2/n3,n3/n4, n4/n6, n3/n5, n5/n7}
    \draw (\from) -- (\to);
\end{tikzpicture}}
\end{center}
\caption{\small{A bipartite graph  and its  2-colouring.}} \label{fig: bipartite graph}
\end{figure}
\end{example}

\begin{remark}
\rm For any cell graph it follows from  the description above that  at any critical point the Hessian of an admissible ${\bf S}^1$-invariant function $h$ is a weighted Laplacian. When $\delta$ is an even function, the weight 
at  each $(i,j)$-edge is $\delta^{\prime \prime}(\theta_i - \theta_j).$ If $\delta$ is not even,  then $\G$ is bipartite  and $h$ is given by (\ref{eq: bipartite S1-inv}); in this case the weight at the $(i,j)$-edge is still  $\delta^{\prime \prime}(\theta_i - \theta_j)$, but now  $i \in \V_k$. 
  \end{remark}

It is well known that without symmetry constraints, a smooth function on a manifold generically  is a Morse function.  In the present case, the Hessian, being a weighted Laplacian,  has a  zero eigenvalue (with eigenvector with all elements equal to 1).  We then have: 

\begin{proposition} \label{prop: degenerate critical points}
For any  cell graph $\G$, all critical points  of $\S^1$-invariant admissible functions on $\G$ are degenerate,  so no Morse singularities are to be expected.  
\end{proposition}

Let us now remark that simultaneous existence of  totally synchronous and 2-colour critical points 
can be expected due to extra symmetries; see Remark~\ref{rmk:coexistence}.  In view of Proposition~\ref{prop: generic 2-colour}, in the search of  2-colour critical points bipartiteness of $\G$ is a generic assumption. In this case, $\V = \V_1 \mathbin{\dot{\cup}} \V_2$ and  2-colour critical points are given by 
$\bar{\theta} = (\theta_1, \ldots, \theta_n)$, 
\[ \theta_i \ = \   \left\{ \begin{array}{ll} 
\alpha_0,  & i \in \V_1  \\
\beta_0, & i \in \V_2, \end{array} \right. \]
$\alpha_0 \neq \beta_0$. Then, we have  synchrony and a 2-colour configuration simultaneously  as long as the coupling function has two distinct critical points, one of them being the origin. In fact, $\theta^0 = (\theta_0, \ldots, \theta_0)$ is critical for $h$ if, and only if,  ${\delta}^\prime(0) = 0$ and, from (\ref{eq: S1 critical bipartite}), $\bar{\theta}$ is critical for $h$ if, and only if, $\delta^\prime(\alpha_0 - \beta_0) = 0.$ In addition, the Hessians at $\theta^ 0$ and at $\bar{\theta}$ are the weighed Laplacians with weights $\delta^{\prime \prime}(0)$ and $\delta^{\prime \prime}(\alpha_0 - \beta_0)$; that is, if ${\cal L}$ is the graph Laplacian of $\G$, then 
\[ H(\theta^0) \ = \ \delta^{\prime \prime}(0) {\cal L}, \ \ \  H(\bar{\theta}) \ = \ \delta^{\prime \prime}(\alpha_0- \beta_0) {\cal L}, \]
and hence the nature of one pair is determined by the nature of the other. The aim of the next subsection is to consider the case when $\G$ is a ring and the coupling function is under some special constraints.  We shall encounter further critical points co-existing with these two, and their nature shall also be investigated. 

\subsection{Stable synchronous patterns on a ring of many cells} \label{subseq: minimum on a ring}

We assume that the graph is a ring of $n$ cells, so the  admissible function is given by (\ref{eq: function on a ring}):
\[ h(\theta) \ = \ \sum_{i=1}^{n} \delta(\theta_i -  \theta_{i-1}), \]
with $\theta_{0} = \theta_n$. By rescaling, we consider here ${\bf S}^1 \equiv [0,1]$ and  the domain of $\delta$  in (\ref{eq: S1-invariant beta}) to be the unit interval, $\delta : [0,1] \to \r.$ Assume that $\delta$ can be extended to an even smooth function to $\r$.
On the configuration space $\r^n$, the gradient of $h$ is an admissible vector field on the $n$-ring. In applications, $h$ represents the total energy associated to the gradient system
\begin{equation} \label{eq: gradient system}
 \dot{\theta} \ = \ - \nabla h(\theta), 
 \end{equation}
whose stable equilibria correspond to minimum points of $h$. Based upon that, the global minimum of $h$ shall be called here the ground state. From Proposition~\ref{prop: degenerate critical points}, the critical points of $h$ are not Morse singularities, for its Hessian at any point  is a weighted Laplacian, and so  the kernel is at least one dimensional. With that in mind, we shall say that an equilibrium 
of (\ref{eq: gradient system}) is {\it asymptotically stable}, or simply {\it stable},  if all nonzero eigenvalues of the Hessian of $h$ at this point are positive.  \\

As mentioned previously, in many applications the coupling is the cosine function; see  \cite{Bronski et al,Korshunov, Lee, MS, Walter Chatelain}. Here we consider a slightly more general class of coupling,  but  chosen suitably to provide a clear understanding of  how it combines with the architecture of the graph to reveal the ground state for the system other than the obvious alternating up-down configuration, which corresponds to antiferromagnetism interactions mentioned earlier.  
More precisely, we  assume that the coupling function $\delta$ satisfies the following conditions: 

(C1) $\delta$ is even and 1-periodic; 

(C2) The first derivative $\delta^{\prime}$ vanishes only at 0 and 1/2; 

(C3) The second derivative $\delta^{\prime \prime}$  is monotone 
in the interval $(0,1/2)$. \\

We now introduce the variables
\[ x_i \ = \ \theta_{i+1} - \theta_i, \ \ 1 \leq i \leq n, \]
where $\theta_{n+1} = \theta_1$.  From (C1)  we can restrict our analysis to  $ x_i \in [0,1/2]$, for $i=1, \ldots, n.$   Critical points of $h$ are then given by 
$\delta^{\prime} (x_{i+1}) \ = \  \delta^{\prime} (x_{i}),   1 \leq i \leq n-1  $, see (\ref{eq: cr pts on a ring}),
where
\[  \sum_{i=1}^{n} x_i \in \z. \]
For $m$ denoting the sum above it follows that  $ 0 \leq m \leq n/2$. In addition, due to the $\S^1$-invariance of $h$ it suffices to describe its critical points  on ${\bf T}^n$.  In fact, for the representation of the dihedral group ${\bf D}_n$ of permutations that preserve the $n$-ring, this is a 
${\bf D}_n \times \S^1$-invariant function, for the  ${\bf D}_n \times \S^1$-action on ${\bf T}^n$ given by
\[ (\pi, \varphi) \cdot (\theta_1, \ldots, \theta_n) \ = \  (\theta_{\pi^{-1}(1)} + \varphi, \ldots, \theta_{\pi^{-1}(n)} + \varphi). \]
In Table~\ref{table critical points} we present the list of all possible equilibria up to their isotropies in  ${\bf D}_n \times \S^1$. We denote by  ${\bf Z}_n(\tau, \varphi)$ the isotropy subgroup generated by $(\tau, \varphi)$, where   $\tau$ is the permutation 
\[ \tau \ = \  \left( \begin{array}{cccc} 
1 & 2 & \cdots & n \\
2 & 3 & \cdots &  1 \\ \end{array} \right) \]
and $\varphi = m/n$, $m=0, \ldots, \lfloor n/2 \rfloor$. 
If we label the cells clockwise,  a ${\bf Z}_n(\tau, \varphi)$-isotropy equilibrium is a configuration such that each cell has a phase shift by   $\varphi$ with respect to its prior neighbour (Figure~\ref{fig: ring} on the left), whereas at a trivial isotropy equilibrium each cell has a phase shift by $\xi$ or $\eta$ relatively to its prior neighbour, where $\xi \neq \eta \in  [0,1/2]$. As we shall see in Proposition~\ref{prop: stable equilibria} below,  a trivial isotropy equilibrium  is stable only if  its configuration is of the form given in Figure~\ref{fig: ring} on the right for $\xi > \eta$.

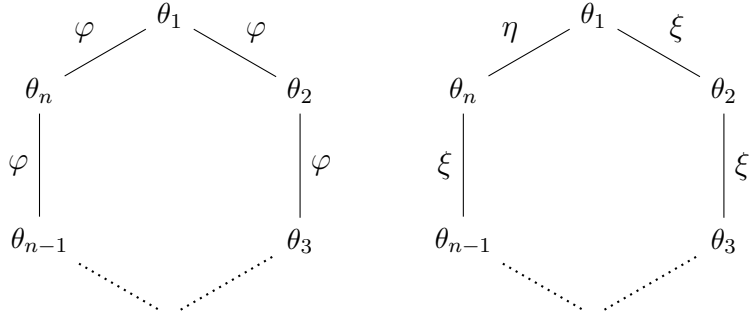
\begin{figure}[H]  
\begin{center}
\begin{tikzpicture}[auto, bend angle=60]
\node (a) at (90:2) {\small{$\theta_1$}}; 
\node (b) at (30:2) {\small{$\theta_2$}};
\node (c) at (150:2) {\small{$\theta_n$}};
\node (d) at (210:2) {\small{$\theta_{n-1}$}};
\node (e) at (270:2) {\small{}};
\node (f) at (330:2) {\small{$\theta_3$}};
 \draw (b) to node {} node [swap] {$\varphi$}  (a);
 \draw (f) to node {} node [swap] {$\varphi$}  (b);
 \draw [dotted, thick] (e) to node {} node [swap] {}  (f);
\draw [dotted, thick] (d) to node {} node [swap] {} (e);
\draw (c) to node {} node [swap] {$\varphi$}  (d);
\draw (a) to node {} node [swap] {$\varphi$}  (c);
\end{tikzpicture} \ \ \ \ \ \ \ 
\begin{tikzpicture}[auto, bend angle=60]
\node (a) at (90:2) {\small{$\theta_1$}}; 
\node (b) at (30:2) {\small{$\theta_2$}};
\node (c) at (150:2) {\small{$\theta_n$}};
\node (d) at (210:2) {\small{$\theta_{n-1}$}};
\node (e) at (270:2) {\small{}};
\node (f) at (330:2) {\small{$\theta_3$}};
 \draw (b) to node {} node [swap] {$\xi$}  (a);
 \draw (f) to node {} node [swap] {$\xi$}  (b);
 \draw [dotted, thick] (e) to node {} node [swap] {}  (f);
\draw [dotted, thick] (d) to node {} node [swap] {} (e);
\draw (c) to node {} node [swap] {$\xi$}  (d);
\draw (a) to node {} node [swap] {$\eta$}  (c);
\end{tikzpicture}
\end{center}
\caption{\small{On the left, a ${\bf Z}_n(\tau, \varphi)$-isotropy equilibrium, for which $\theta_{i+1} = \theta_{i} + \varphi$, $i=1, \ldots, n$ with $\theta_{n+1} = \theta_1$. On the right, a ${\bf 1}$-isotropy equilibrium, for which  
$\theta_{i+1} = \theta_{i} + \xi$, $i=1, \ldots, n-1$ and   $\theta_{1} = \theta_n + \eta$, $\eta \neq \xi.$}} \label{fig: ring}
\end{figure}
It is  a direct consequence of the conditions (C1)-(C3) that there are no equilibria rather than those in Table~\ref{table critical points}.  Also, we notice that trivial isotropy equilibrium points only occur for special couplings. In fact, the values $\xi$ and $\eta$ in the last row of this table are such that
\begin{equation} \label{eq: x and y}
 0 \leq \eta < \xi \leq 1/2, \ \ \delta^\prime (\xi) = \delta^\prime(\eta), 
 \end{equation}
whose existence is guaranteed by the conditions imposed on $\delta$. However,  these will correspond to equilibria of the system  
if and only if $p \xi + q \eta$ is an integer $m$, for $p,q$ and $m$ in the conditions presented therein. This is the case for example for  the values $\xi = 1/2$, $\eta=0$, giving rise to equilibrium points  $\theta = (\theta_1, \ldots, \theta_n)$ such that there is en even number of pairs of neighbours whose states differ by 1/2.

\begin{table}[H] 
{\small
\centering
\begin{tabular}{|c|c|c|c|c|}
\hline\hline
$m \in \z$,& Isotropy & Critical point & Value of & Laplacian \\ 
$0 \leq m \leq  n/2 $ & &representative & energy  & weights \\ \hline
0 & ${\bf D}_n$ & $(0, 0, \ldots, 0)$ & $n \ \delta(0)$ & $\delta^{\prime \prime}(0)$ \\ \hline
$ \neq 0, n/2$ & ${\bf Z}_n(\tau, m/n)$ & $(0, m/n, 2m/n, \ldots, m)$ & $n \ \delta(m/n)$ & $\delta^{\prime \prime}(m/n)$\\ \hline
$n/2$ (if $n$ even) &   ${\bf Z}_n(\tau, 1/2)$  & $(0, 1/2, 0, \ldots, 1/2)$ & $n \ \delta(1/2)$ & $\delta^{\prime \prime}(1/2)$\\ \hline
$ p\xi+q\eta$, & {\bf 1} &$(\theta_1, \ldots, \theta_n)$  &  $p \delta(\xi) + q \delta (\eta)$ & $\delta^{\prime \prime}(\xi),$ \\
for  $p, q \in {\n},$   $p+q =n$ &   & $\theta_{i+1} - \theta_{i} = \xi$ or $\eta$  & &  $\delta^{\prime \prime}(\eta)$ \\ \hline \hline 
\end{tabular}
\caption{\small{Critical points of admissible functions on a ring of $n$ cells.}} \label{table critical points}}  
\end{table}

By direct investigation, if   $\delta^{\prime \prime}(0) >0$ then the fully synchronous pattern  $x_i = 0$, $1 \leq i \leq n$, is the ground state. This corresponds to  the ferromagnetic state and it falls into the case $\alpha + \beta = 0, \beta < 0$   discussed in Subsection~\ref{subsec: tot synchr regular}. \\

 On the other hand, if  $\delta^{\prime \prime}(0) <0$ then fully synchrony corresponds to a totally unstable equilibrium. The aim here is to extract from Table~\ref{table critical points} all  possible stable equilibria in this case. In addition, we detect the ground state if we allow the ring to have a large number of cells.  In this context, we shall then consider:  \\
 
(C4)  $\delta^{\prime \prime}(0) < 0$ and $\delta^{\prime \prime}(1/2) >0$. \\

The results are given in Propositions~\ref{prop: stable equilibria} and \ref{prop: ground state} and rely on the fact that the Hessian of admissible functions defined on rings is, at any point, a graph weighted Laplacian. 
We shall see how stability is attained depending  on the topoloy of the graph combined with its weights. The weights are real numbers associated with each edge, so as expected they are determined by  the constrains imposed on the coupling function $\delta$.

We first recall an important fact of \cite{Bronski DeVille} about eigenvalues of weighted Laplacians.  For any weighted graph $\G$ with $n$ vertices, let $\G_+$ (resp. $\G_-$) denote  the subgraph  with the same vertex set as $\G$ together with the edges of positive (resp. negative) weights. We also define the three indices $n_0, n_-$ and $n_+$ as the numbers of zero, negative and positive eigenvalues of the weighted Laplacian.  In \cite{Bronski DeVille} the authors give the best possible bounds on these three indices  involving only topological information, i.e., connectivity of the graph and the sign information on the edge weights. For $c(\G)$ denoting the number of connected component of any graph $\G$, the estimates are given as follows:
\begin{equation} \label{eq: bounds}
\begin{tabular}{rl}
$c(\G_+) - 1 \leq$  & $n_-$  $\leq n - c(\G_-)$ \\
$c(\G_-) - 1 \leq $ & $ n_+$  $\leq n - c(\G_+)$ \\  
$ 1 \leq$  & $ n_0 \ \leq  n+2 - c(\G_-) - c(\G_+)$.
\end{tabular}
\end{equation}

\begin{proposition} \label{prop: stable equilibria}
Let $\G$ be a ring of $n$ cells coupled under conditions (C1)-(C4). Then the stability of critical points that appear in Table~\ref{table critical points} is as follows: \\

\noindent {\rm (1)} For $n$ even, an equilibrium with ${\Z}_n(\tau, 1/2)$-isotropy is stable;

\noindent {\rm (2)} For $m \neq 0, n/2$, an equilibrium with ${\Z}_n(\tau, m/n)$-isotropy is stable if and only if $\delta^{\prime \prime}(m/n) >0$;

\noindent {\rm (3)} For $n$ odd, a trivial isotropy equilibrium $\theta = (\theta_1, \ldots, \theta_n)$, $\theta_i \in \{0,1/2\}$, is stable only if at most one pair of neighbours assume the same state value. 

\noindent {\rm (4)} Possible trivial isotropy equilibria with $\xi, \eta$ in the open interval $(0,1/2)$ are stable only if $p=n-1$, $q=1$.
\end{proposition}

\dem The weighted Laplacian at  the ${\Z}_n(\tau, 1/2)$-isotropy equilibrium (for $n$ even) is  given by
$\delta^{\prime \prime}(1/2) {\cal L}$, where ${\cal L}$ is the graph Laplacian (without weights). The $n-1$  nonzero eigenvalues of $L$ are positive, so stability of this point follows from (C4).    The stability of 
${\Z}_n(\tau, m/n)$-isotropy points follows analogously.  Consider now the trivial isotropy points 
$\theta$ given in (3).  The weights in this case are given by $\delta^{\prime \prime}(0)$ and 
$\delta^{\prime \prime}(1/2)$. Now, from (\ref{eq: bounds})  a necessary condition for any point to be stable is $c(\G_+) = 1$, i.e., $\G_+$ must be connected. Under condition (C4) this is the case on a ring  if and only if  the states of cells alternate between $0$ and $1/2$ except for possibly one pair of neighbours.  Finally we consider $\theta$ given in (4). We have that 
$\delta^{\prime \prime}(\xi) >0$ and $\delta^{\prime \prime}(\eta) < 0$ by (C3), and  the result follows analogously as for case (3).
\hfill \eproof

\begin{remark} \normalfont  \label{rmk: frustration}
Case (3) of Proposition~\ref{prop: stable equilibria} corresponds to the {\it frustration} phenomenon in magnetic systems. Geometrical frustration is a feature in magnetism related to the arrangement of spins. In a ring with an even number of spins, the energy is minimized when each spin is aligned opposite to neighbours (see Proposition~\ref{prop: ground state} below). If the number of cells is three for example, once the first two spins align anti-parallel, the third one is `frustrated' because of its two possible orientations, up and down, which give the same energy. The third spin cannot simultaneously minimize its interactions with both of the other two. 
From Proposition~\ref{prop: stable equilibria} (3), this is the case if there is exactly one pair of neighbours with same state value, as illustrated in Figure~\ref{fig: frustration}.
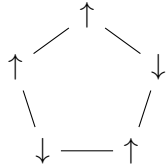
\begin{figure}[H] 
\begin{center}
\begin{tikzpicture}[auto, bend angle=72]
\node (a) at (90:1) {\small{{\color{black} $\uparrow$}}}; 
\node (b) at (18:1) {\small{{\color{black} $\downarrow$}}};
\node (c) at (162:1) {\small{{\color{black} $\uparrow$}}}; 
\node (d) at (234:1) {\small{{\color{black} $\downarrow$}}}; 
\node (e) at (306:1) {\small{{\color{black} $\uparrow$}}}; 
 \draw (b) to node {} node [swap] {}  (a);
\draw (a) to node {} node [swap] {}  (c);
\draw (c) to node {} node [swap] {}  (d);
\draw (d) to node {} node [swap] {}  (e);
 \draw (b) to node {} node [swap] {}  (e);
\end{tikzpicture} 
\end{center}
\caption{\small{The frustration phenomenon in a ring of five cells.}} \label{fig: frustration}
\end{figure}
\end{remark}

\begin{proposition} \label{prop: ground state}
Let $\G$ be a ring of $n$ cells coupled under conditions (C1)-(C4).  If $n$ is even then the 
ground state $(\theta_1, \ldots, \theta_n)$ of (\ref{eq: gradient system})  is  given by $\theta_{i+1} - \theta_i = 1/2,$ $i=1, \ldots, n$. If $n$ is odd and sufficiently large, then the ground state is  given by $\theta_{i+1} - \theta_i = 1/2 - 1/2n$, $i=1, \ldots, n.$ 
\end{proposition}

\dem For $n$ even, the conclusion follows by direct investigation on the values of the energy function $h$ given in Table~\ref{table critical points}. For $n$ odd we have to compare  among ${\Z}_n(\tau, m/n)$-isotropy
and trivial isotropy stable equilibria given by Proposition~\ref{prop: stable equilibria}.  Since $\delta$ is a decreasing function on $[0,1/2]$, the lowest energy with ${\Z}_n(\tau, m/n)$ isotropry is attained for the largest $m$, namely $m = (n-1)/2$. This corresponds to the equilibrium 
$\tilde{\theta} = (\theta_1, \ldots, \theta_n)$,  $\theta_i = 1/2 - 1/2n$, $i=1, \ldots, n.$  Hence, it remains to compare the energy values at this point and at trivial isotropy points $\bar{\theta}$ which are given by (3) and (4) in Proposition~\ref{prop: stable equilibria}.  

We have that $(n-1) \xi + \eta = m,$ $0 \leq m \leq n/2$ and $\eta = \xi - \epsilon$, $0 < \epsilon < 1/2$.   For $n$ odd, $m \leq (n-1)/2$. If $m < (n-1)/2$, then
\[ (n-1)\xi + \eta \leq \frac{n-1}{2} -1, \]
so 
\[ \xi \leq  \frac{n-1}{2n} + \frac{\epsilon -1}{n} < \frac{1}{2} - \frac{1}{2n}. \]
Hence, $\delta(\xi) > \delta(1/2 - 1/(2n)),$ and since $\delta(\eta) > 0$,  
\[ n  \ \delta\bigl(\frac{1}{2} - \frac{1}{2n}\bigr) < (n-1) \delta(\xi) + \delta(\eta).\]
If $m = (n-1)/2,$ then 
\[ \xi  = \frac{n-1}{2n} + \frac{\epsilon}{n}, \  \ \eta = \frac{n-1}{2n} + \frac{1-n}{n} \epsilon. \] 

If  $n$ is sufficiently large, then
 \[ \delta\bigl(\frac{1}{2} - \frac{1}{2n}\bigr) < 0, \ \ 
  \delta \bigl( \frac{1}{2} - \frac{1}{2n}  + \frac{\epsilon}{n}\bigr) \approx \delta\bigl( \frac{1}{2} - \frac{1}{2n} \bigr). \]
We now just compare
\[ f(\tilde{\theta}) \ = \ n \delta\bigl( \frac{1}{2} - \frac{1}{2n} \bigr) = (n-1)\delta\bigl( \frac{1}{2} - \frac{1}{2n} \bigr) + \delta\bigl( \frac{1}{2} - \frac{1}{2n} \bigr), \]
with
\[ f(\bar{\theta}) \ = \  (n-1) \delta(\xi) + \delta(\eta) \]
to conclude that, in this limit, $h(\tilde{\theta}) < h(\bar{\theta}).$ 
 \hfill \eproof  

Proposition~\ref{prop: ground state} asserts that for the bipartite case the ground states are achieved in the up-down configuration, see Figure~\ref{fig: ground states}(a). This is the corresponding result for rings to the analysis in  \cite{Lee}  of the classical AFXY model, for which bipartite lattices of spins attain the ground state when the two sublattices are aligned in opposite directions. For the non bipartite case,  (when $n$ is odd), the ground state value depends on $n$  and tends to the up-down arrangement in the limit $n \to \infty$, as depicted in  Figure~\ref{fig: ground states}(b).
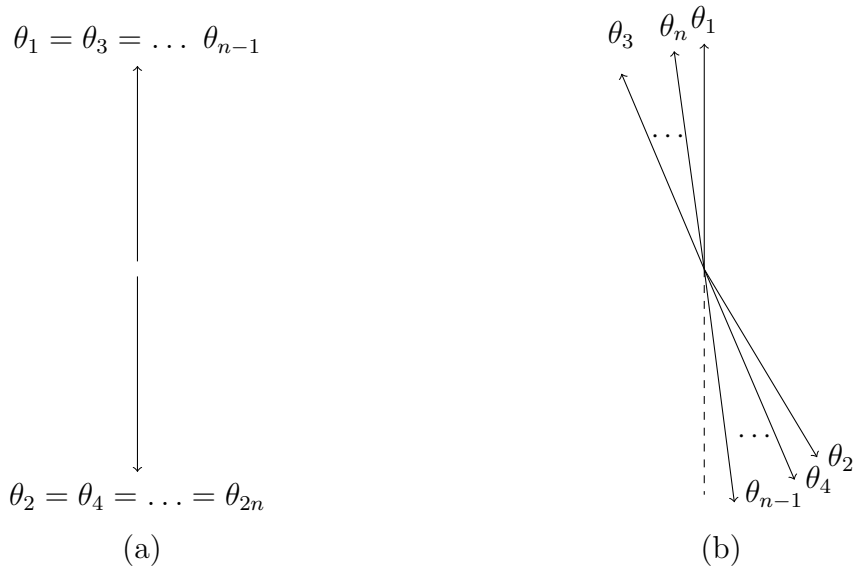
\begin{figure}[H] 
\begin{center}
\begin{tikzpicture} 
\foreach \coordinate/\label/\pos in {{(0,2.7)/$\theta_1=\theta_3=\ldots \ \theta_{n-1}$/above},{(0,-2.7)/$\theta_2=\theta_4=\ldots=\theta_{2n}$/below}}
\node[\pos] at \coordinate {\label};\draw[<-] (0,-2.7) -- (0,-0.1);
\draw[->] (0,0.1) -- (0,2.7);
\end{tikzpicture} 
\hspace*{4cm}
\begin{tikzpicture} 
\draw[black,->] (0,0) -- (0,3);
\draw[black,->] (0,0) -- (-1.1,2.6);
\draw[black,->] (0,0) -- (1.5,-2.5);
\draw[black,->] (0,0) -- (1.2,-2.8);
\draw[black,->] (0,0) -- (0.4,-3.1);
\draw[black,->] (0,0) -- (-0.4,2.9);
\draw [dashed] (0, 0) -- (0,-3);
\foreach \coordinate/\label/\pos in {{(0,3)/$\theta_1$/above},{(-1.1,2.8)/$\theta_3$/above},{(1.5,-2.5)/$\theta_2$/right},{(1.2,-2.8)/$\theta_4$/right},{(0.4,-3.0)/$\theta_{n-1}$/right},{(0.7,-2)/{\small $\cdots$}/below},{(-0.45,2)/{\small $\cdots$}/below},{(-0.4,2.9)/$\theta_n$/above}}
\node[\pos] at \coordinate {\label};
\end{tikzpicture} 

\hspace*{4cm} (a) \hfill (b) \hspace*{4cm}
\end{center}
\caption{\small{Ground state in a ring of $n$ cells under conditions (C1)-(C4) when (a) $n$ is even and (b) $n$ is odd and large.}} \label{fig: ground states}
\end{figure}

\vspace*{1cm}

\noindent {\it Acknowledgments.} The research of M.M. was supported by FAPESP,  BPE grant 2013/11108-7.

\end{document}